\documentclass{article}
\usepackage{amsmath,amsfonts,amsthm}

\newcommand{\N}{\mathbb{N}}
\newcommand{\Z}{\mathbb{Z}}
\newcommand{\R}{\mathbb{R}}
\newcommand{\C}{\mathbb{C}}
\newcommand{\B}{\mathbb{B}}
\newcommand{\coneS}{\mathcal{S}}
\newcommand{\sbody}[2]{\textstyle \frac{#1}{#2} }
\newcommand{\rmref}[1]{(\romannumeral 0\ref{#1}) 
 }

\newtheorem{Th}{Theorem}[section]
\newtheorem{Lem}[Th]{Lemma}
\newtheorem{corollary}[Th]{Corollary}
\newtheorem{definition}[Th]{Definition}
\newtheorem{remark}[Th]{Remark}

\newenvironment{Proof}[1][Proof.]{\begin{trivlist}
\item[ \hskip \labelsep {\bfseries #1}]} {\end{trivlist}}

\begin{document}

\title{Stieltjes, Poisson and other integral representations for \\
functions of Lambert $W$}
\author{German A. Kalugin, David J. Jeffrey \& Robert M. Corless \\
Department of Applied Mathematics, The University of Western Ontario,\\
London, Ontario, Canada}
\maketitle

\begin{abstract}
We show that many functions containing $W$ are Stieltjes functions.
Explicit Stieltjes integrals are given for functions $1/W(z)$, $W(z)/z$, and others.
We also prove a generalization of a conjecture of Jackson, Procacci \& Sokal.
Integral representations of $W$ and related functions are also given which are associated with
the properties of their being Pick or Bernstein functions.
Representations based on Poisson and Burniston--Siewert integrals are given as well.
\end{abstract}

\section{Introduction}
The Lambert $W$ function is the multivalued inverse of the mapping $W \mapsto W e^W$.
The branches, denoted by $W_k$ ($k\in\Z$), are defined through the equations \cite{Big paper}
\begin{eqnarray} \label{eq:def}
\forall z\in \C,\qquad W_k(z)\exp(W_k(z)) &=& z\ , \\
\label{eq:branch}
W_k(z) &\to & \ln_k z  \mbox{ as } \Re z \to \infty\ ,
\end{eqnarray}
where $\ln_k z=\ln z + 2\pi i k$, and $\ln z$ is the principal branch of
natural logarithm~\cite{Unwinding}.
This paper considers only the principal branch $k=0$, which is the branch that maps the positive real
axis onto itself, and therefore we shall usually abbreviate $W_0$ as $W$ herein.

The $W$ function is rich in integral representations. As a consequence, $W$ is a member of a
number of function classes, specifically, the classes of Stieltjes functions, Pick functions and
Bernstein functions, including subclasses Thorin-Bernstein functions and complete Bernstein functions.
A description of the classes can be found in a review paper \cite{Berg2008} and a recently published
book \cite{Schilling}.
The integral representations are related to the fact that $W$ is a real symmetric function,
in the terminology of \cite[p.\,160]{Baker} (see also \cite[p.\,155]{Titchmarsh}),
with positive values on the positive real line.
In addition to integral representations connected to the above function classes,
integrals following Poisson~\cite{Poisson} and Burniston--Siewert~\cite{Burniston} are given.

The classes of Stieltjes functions and Bernstein functions are tightly connected
with the class of completely monotonic functions that have many applications in different
fields of science;
the list of appropriate references is given in \cite{AlzerBerg1}.
Therefore we shall also study the complete monotonicity of some functions containing $W$.

The properties and integral representations mentioned above have interesting computational implications.
For example,
that $W(z)/z$ is a Stieltjes function means that the poles of successive Pad\'e approximants interlace
and all lie on the negative real axis \cite[p.\:186]{Baker} (here in the interval $-\infty<z<-1/e$).
In addition, some of the integral representations permit spectrally convergent
quadratures for numerical evaluation.

\subsection{Properties of $W$}
For convenience, we recall from~\cite{Big paper} some properties of $W$ that are used below.
The function is discontinuous on the branch cut $\B \subset \R$, defined to be the interval
$\B=(-\infty, -1/e]$. On the cut plane $\C\backslash \B$, the function is
holomorphic.
Its real values obey $-1\le W(x) <0$ for $x\in[-1/e,0)$, $W(0)=0$ and $W(x)>0$ for $x>0$.
The imaginary part of $W(t)$ has the following range of values for real $t$.
\begin{equation} \label{eq:ImW}
 \Im W(t)\in(0,\pi) \mbox{   for   } t\in(-\infty,-1/e) \mbox{   and   } \Im W(t)=0 \mbox{   otherwise.}
\end{equation}
$\Im W(z)$ and $\Im z$ have the same sign in the cut plane $\C\backslash\R$, or equivalently
\begin{equation} \label{signImW}
 \Im W(z)\Im z >0 \ .
\end{equation}
$W$ has near conjugate symmetry, meaning $W(\overline z) = \overline {W(z)}$ except
on the branch cut $\B$.
We also note
\begin{equation} \label{eq:Wright}
 W(z) = \ln z - \ln W(z)
\end{equation}
in the cut plane $\C\backslash(-\infty,0]$.
The Taylor series near $z=0$ is
\begin{equation} \label{eq:Taylor}
 W(z) =\sum_{n=1}^\infty(-n)^{n-1}\frac{z^n}{n!}
\end{equation}
with radius of convergence $1/e$, while the
asymptotic behaviour of $W(z)$ near its branch point is given by
\begin{equation} \label{W:z=-1/e}
W(z)\sim -1+\sqrt{2(ez+1)} \mbox{\quad as \quad} z\rightarrow -1/e \ .
\end{equation}
It follows from \eqref{eq:Taylor} and \eqref{eq:branch} that
\begin{equation} \label{WnearZero}
 W(z)/z \rightarrow 1 \quad \mbox{as} \quad z\rightarrow 0 \ .
\end{equation}
\begin{equation} \label{WnearInfty}
 W(z)/z \rightarrow 0 \quad \mbox{as} \quad z\rightarrow \infty \ .
\end{equation}
If $z=t+is$ and $W(z)=u+iv$, then
\begin{equation*}
e^u(u\cos v-v\sin v)=t , \quad  e^u(u\sin v+v\cos v)=s \ .
\end{equation*}
For the case of real $z$, i.e. $s=0$, the functions $u=u(t)$ and $v=v(t)$ are defined by
\begin{align}
u&=-v\cot v ,          \label{eq:u(v)} \\
t&=t(v)=-v\csc(v)e^{-v\cot v} \label{eq:t(v)}  \ .
\end{align}
For the case of purely imaginary $z$, i.e. $t=0$, the functions $u=u(s)$ and $v=v(s)$ obey
\begin{align}
u&=v\tan v,             \label{Imz:u(v)} \\
s&=s(v)=v\sec(v)e^{v\tan v}  \label{Imz:s(v)} \ .
\end{align}
The derivative of $W(z)$ is given by
\begin{equation} \label{eq:DW}
 W^\prime(z)=\frac{W(z)}{z(1+W(z))} \ .
\end{equation}
It follows from \eqref{WnearZero} and \eqref{WnearInfty}
\begin{equation} \label{DWnearZero}
 W^\prime(z)\rightarrow 1 \quad \mbox{as} \quad z\rightarrow 0 \ ,
\end{equation}
\begin{equation} \label{DWnearInfty}
 W^\prime(z)\rightarrow 0 \quad \mbox{as} \quad z\rightarrow\infty \ .
\end{equation}
Near conjugate symmetry implies
\begin{equation} \label{symDW}
dW(\bar{z})/d\bar{z}=\overline{dW(z)/dz} \
\end{equation}
for $z\notin(-\infty,0]$.
%

\subsection{Stieltjes functions} \label{sec:properties}
We now review the properties of Stieltjes functions, again concentrating on results that will
be used in this paper.
We must note at once that
there exist several different definitions of Stieltjes functions in the literature, and
here we follow the definition of Berg~\cite{Berg2008}.
\begin{definition} \label{def:Stieltjes fun}
 A function $f:(0,\infty)\rightarrow\R$ is called a
\emph{Stieltjes function} if it admits a representation
\begin{equation} \label{def:realSt}
f(x)=a+\int_0^\infty\frac{d\mu(t)}{x+t} \mbox{\quad for \quad } x>0\ ,
\end{equation}
where $a$ is a non-negative constant and $\mu$ is a positive measure on $[0,\infty)$
such that $\int_0^\infty(1+t)^{-1}d\mu(t)<\infty$.
\end{definition}

A Stieltjes function is also called a \emph{Stieltjes transform} \cite[p.\,127]{BergForst}, and
the integral in the right-hand side of \eqref{def:realSt} is called a Stieltjes integral of $f(x)$.
Except in \S \ref{sec:W-S} below, the term Stieltjes function will here always refer to
Definition \ref{def:realSt}.

\begin{Th}
 \label{th:StieltjesProp}

The set $\coneS$ of all Stieltjes functions forms a
convex cone \cite[p.\,127]{BergForst} and possesses the following properties.

\begin{enumerate}
\renewcommand{\labelenumi}{(\roman{enumi})}
\item $f\in\coneS\:\backslash\left\{0\right\}\Rightarrow\frac{1}{f(1/x)}\in\coneS$ \label{en:1}
\item $f\in\coneS\:\backslash\left\{0\right\}\Rightarrow\frac{1}{xf(x)}\in\coneS$ \label{en:2}
\item $f\in\coneS\Rightarrow\frac{f}{cf + 1}\in\coneS \quad (c\geq 0)$ \label{en:3}
\item $f,g\in\coneS\:\backslash\left\{0\right\}\Rightarrow f\circ\frac{1}{g}\in\coneS$ \label{en:4}
\item $f,g\in\coneS\:\backslash\left\{0\right\}\Rightarrow \frac{1}{f\circ g}\in\coneS$ \label{en:5}
\item $f,g\in\coneS\Rightarrow f^\alpha g^{1-\alpha}\in\coneS \quad (0\leq\alpha\leq1)$ \label{en:6}
\item $f\in\coneS\Rightarrow f^\alpha\in\coneS \quad (0\leq\alpha\leq1)$ \label{en:7}
\item $f\in\coneS\:\backslash\left\{0\right\}\Rightarrow\frac{1}{x}\left(\frac{f(0)}{f(x)}-1\right)\in\coneS$ \label{en:8}
\item $f\in\coneS\:\backslash\left\{0\right\}, \lim_{x\rightarrow0}xf(x)=c>0 \Rightarrow f(x)-c/x\in\coneS$ \label{en:9}
\item $f\in\coneS\Rightarrow f^\alpha(0)-f^\alpha(1/x)\in\coneS \quad (0\leq\alpha\leq1)$ \label{en:10}
\item $f\in\coneS\:\backslash\left\{0\right\}\Rightarrow\frac{1}{x}\left(1-\frac{f(x)}{f(0)}\right)\in\coneS$ \label{en:11}
\item $f\in\coneS, \lim_{x\rightarrow\infty}f(x)=c>0\Rightarrow(c^\beta-f^\beta)\in\coneS \quad (-1\leq\beta\leq0)$ \label{en:12}
\end{enumerate}
\end{Th}

\begin{Proof}
Properties \rmref{en:1}--\rmref{en:7} are listed in \cite {Berg2008};
property \rmref{en:6} is due to the fact that
the Stieltjes cone is logarithmically convex~\cite{Berg1979}
and property \rmref{en:7} is its immediate consequence.
Property \rmref{en:8} is taken from  \cite[p.\,406]{BenderOrszag}.
Property \rmref{en:9} follows from properties \rmref{en:2} and \rmref{en:8} in the following way:
$f\in\coneS\:\backslash\left\{0\right\} \Rightarrow g(x)=1/(xf(x))\in\coneS \Rightarrow
(g(0)/g(x)-1)/x=(xf(x)/c-1)/x\in\coneS \Rightarrow f(x)-c/x\in\coneS$.
The last three properties \rmref{en:10}--\rmref{en:12} will be proved in \S~\ref{sec:cbf}.
\end{Proof}

A Stieltjes function $f$ has a holomorphic extension to the
cut plane $\C\backslash(-\infty,0]$ satisfying $f(\bar{z})=\overline{f(z)}$
(see \cite{Berg1979}, \cite{AlzerBerg} and \cite[p.\,11-12]{Schilling})
\begin{equation} \label{def:complexSt}
f(z)=a+\int_0^\infty\frac{d\mu(t)}{z+t} \quad (\left|\arg z\right|<\pi)\ .
\end{equation}
In addition, a Stieltjes function $f(z)$ in the cut plane $\C\backslash(-\infty,0]$
can be represented in the integral form
\begin{equation} \label{eq:DefStieltFun}
f(z) = \int_0^\infty\frac{d\Phi(u)}{1+uz} \quad (\left|\arg z\right|<\pi)\ ,
\end{equation}
with a bounded and non-decreasing function $\Phi(u)$.
The integral \eqref{eq:DefStieltFun} is useful for the analysis of Pad\'e approximations
of Stieltjes functions \cite[Ch.\,5]{Baker}
and equivalent to the representation \eqref{def:complexSt},
which follows from the following observation.

According to properties \rmref{en:1} and \rmref{en:2}, if a function $f\in\coneS$
then $f(1/x)/x\in\coneS$ as well
and hence the latter admits representation \eqref{def:realSt}
\[
\frac{1}{x}f\left(\frac{1}{x}\right)=a+\int_0^\infty\frac{d\rho(t)}{x+t} \ ,
\]
which after replacing $x$ with $1/x$ gives
\begin{equation*} \label{equivalence}
f(x)=a+\int_0^\infty\frac{d\rho(t)}{1+xt} \ ,
\end{equation*}
where the constant $a$ can be included in the integral.
Finally, one considers the holomorphic extension of the last integral to the cut plane
$\C\backslash(-\infty,0]$
similar to obtaining \eqref{def:complexSt}.

There are various kinds of necessary and sufficient conditions implying that a function
$f$ is a Stieltjes function.
Some of them are based on the classical results established by R. Nevanlinna, F. Riesz, and Herglotz.
Here we quote two such theorems, taken from \cite[p.\,93]{Akhiezer} and \cite[Theorem 3.2]{Berg2008}.

\begin{Th} \label{th:criterionA}
A function $g(z)$ admits an integral representation in the upper half-plane in the form
\begin{equation} \label{eq:form2}
g(z) = \int_\R\frac{d\Phi(u)}{u-z} \quad (\Im z>0) \ ,
\end{equation}
with a non-decreasing function $\Phi(u)$ of bounded variation on $\R$
(i.e. $\int_\R d\Phi(u)<\infty$ for smooth $\Phi(u)$),
if and only if $g(z)$ is holomorphic in the upper half-plane and
\begin{equation} \label{eq:conditions2}
\Im g(z)\geq0  \quad\mbox{and}\quad   \sup_{1<y<\infty}\left|yg(iy)\right|<\infty \ .
\end{equation}
\end{Th}

To apply Theorem \ref{th:criterionA} to the integral \eqref{eq:DefStieltFun}
one should set $g(z)=-f(-1/z)/z$ (cf. \cite[(6.12) on p.\,215]{Baker}),
then conditions \eqref{eq:conditions2} read as
\begin{equation} \label{cond:f}
\Im {f(-1/z)/z}\leq0  \quad\mbox{and}\quad  \sup_{1<y<\infty}\left|f(i/y)\right|<\infty \ .
\end{equation}

\begin{Th}\label{th:CriterionB}
A function $f: (0,\infty)\rightarrow\R$ is a Stieltjes function if and only if $f(x)\geq0$ for $x>0$ and
there is a holomorphic extension $f(z)$, $z=x+iy$, to the cut plane $\C\backslash(-\infty,0]$
satisfying
\begin{equation} \label{inequalities}
\Im f(z)\leq0 \mbox{ for } \Im z>0.
\end{equation}
\end{Th}

\begin{remark} \label{rem:Herglotz} The inequalities \eqref{inequalities} alone express a necessary condition for $f$ to be a Stieltjes function.
In the terminology of \cite[p.\,358]{BenderOrszag}, a holomorphic function $f(z)$ is called a \emph{Herglotz function}
if $\Im f>0$ when $\Im z>0$, $\Im f=0$ when $\Im z=0$ and $\Im f<0$ when $\Im z<0$.
Thus, for $f$ to be a Stieltjes function it is necessary that $f$ be an anti-Herglotz function
(cf. \cite[p.\,406]{BenderOrszag}).
\end{remark}

\section{Stieltjes functions containing $W(z)$}
In this section we consider a number of functions containing $W(z)$ and prove that they are Stieltjes
functions. We begin with the function $W(z)/z$.
\subsection{The function $W(z)/z$}
Although we could establish the result conveniently by applying one of the criteria stated in
\S \ref{sec:properties}, we nonetheless first present a direct proof
that is of great importance for further investigations.
Moreover, compared with using the criteria above, the present way allows us to make useful observations
which are given in the remarks following the proof and used in further evidence.
\begin{Th} \label{th:W(z)/z}
The function $W(z)/z$ can be represented in the form of a Stieltjes integral \eqref{eq:DefStieltFun}.
\end{Th}

\begin{Proof} From \eqref{WnearZero}, the function $W(z)/z$ is
single-valued and holomorphic in the same domain as $W(z)$, namely
$D=\left\{ z\in\C \,|\, z\notin \B \right\}$,
and can be represented by the Cauchy's integral formula
\begin{equation} \label{eq:CauchyIntegral}
\frac{W(z)}{z}=\frac{1}{2\pi i}\int_C\frac{W(t)}{t(t-z)}dt \ ,
\end{equation}
where $C$ is the standard `keyhole` contour which consists of a small circle
around the branch point $t=-1/e$ of radius, say $r$, and a large circle around the origin of radius,
say $R$; the circles being connected through the upper and lower edges of the cut along the
negative real axis.
Then for sufficiently small $r$ and large $R$ the interior of the contour $C$ encloses any point in $D$.

Let us consider the integral \eqref{eq:CauchyIntegral} in the limit in which
$r\rightarrow 0$ and $R\rightarrow \infty$.
Using asymptotic estimations  \eqref{W:z=-1/e} and \eqref{eq:branch},
it is easily seen that the contributions of each circle to
the integral \eqref{eq:CauchyIntegral} go to zero.
As a result, in accordance with the assignment of values of $W$ function on the branch cut,
the integral becomes
\[
\frac{W(z)}{z}=\frac{1}{2\pi i}\int_{-\infty}^{-1/e}\frac{W(t)}{t(t-z)}dt +
\frac{1}{2\pi i}\int_{-1/e}^{-\infty}\frac{\overline{W(t)}}{t(t-z)}dt\ ,
\]
which reduces to
\begin{equation} \label{eq:res1}
\frac{W(z)}{z}=\frac{1}{\pi}\int_{-\infty}^{-1/e}\frac{\Im W(t)}{t(t-z)}\ dt\ .
\end{equation}
The change of variable $t=-1/u$ transforms the integral \eqref{eq:res1} to the
form \eqref{eq:DefStieltFun}
\begin{equation} \label{eq:result2}
\frac{W(z)}{z} = \int_0^e\frac{d\Phi(u)}{1+uz} \ ,
\end{equation}
where
\begin{equation} \label{eq:measure}
\Phi(u)=\frac{1}{\pi}\int_0^u \Im W\left(-1/t\right)\ dt \mbox{, } u\in[0,e] \ .
\end{equation}
According to \eqref{eq:ImW}
\begin{equation} \label{ImW:bounds}
0\leq\Im W(-1/t)\leq\pi \mbox{ for } t\in[0,e] \ ,
\end{equation}
hence $\Phi^\prime(u)\geq0$, and thus $\Phi(u)$ is a bounded and non-decreasing function.
\end{Proof}

\begin{remark}  The function $W(z)/z$ is a real symmetric function as
is any Stieltjes function (this immediately follows from Definition 2.2),
which just corresponds to the near conjugate symmetry property.
\end{remark}

\begin{remark} It follows from the proof of the theorem that not only is
$W(z)/z$ a Stieltjes function,
it has a smooth measure $\Phi(u)$ defined by \eqref{eq:measure}.
Moreover, all the moment integrals
\[
\int_0^e t^n\,d\Phi(t) \quad (n=0,1,2,\dots)
\]
exist because so do the integrals
\[
\int_0^e t^n\,\Im W\left(-1/t\right)dt
\]
owing to \eqref{ImW:bounds}.
This remark is useful for justifying the use of Pad\'e approximants for the evaluation of $W(z)$ mentioned in \cite{Big paper}.
\end{remark}

Now we apply Theorem \ref{th:criterionA} to the function $f(z)=W(z)/z$ for which conditions \eqref{cond:f} read as
\begin{equation*}
\Im W(-1/z)\geq0  \quad\mbox{and}\quad  \sup_{1<y<\infty}\left|yW(i/y)\right|<\infty \ .
\end{equation*}
The first condition is satisfied by \eqref{signImW} because $\Im (-1/z)$ and
$\Im z$ are of the same sign.
To verify the second condition we set $W(i/y)=u+iv$ and put $s=1/y$ in
\eqref{Imz:u(v)} and \eqref{Imz:s(v)}.
As a result, since $0<v<\pi/2$ for $y>0$, we obtain
\[
\left|yW(i/y)\right|^2=y^2(u^2+v^2)=y^2 v^2(1+\tan^2 v)=y^2 v^2/\cos^2v=e^{-2v\tan v}\leq1.
\]
To extend the result to the lower half-plane $\Im z<0$ it is enough to take the complex conjugate of
both sides of the representation \eqref{eq:DefStieltFun} and use the near conjugate symmetry of $W$.
Thus Theorem \ref{th:criterionA} gives us one more way to prove that $W(z)/z$ is a Stieltjes function.

\subsection{Other functions}
By Theorem \ref{th:W(z)/z}, $W(x)/x\in\coneS$.
Using this result and the properties of the set $\coneS$ listed in section \ref{sec:properties}
we now give some classes of functions that are members of $\coneS$.

\begin{Th} \label{th:coneS}
The following functions belong to the set $\mathcal{S}$, for $x>0$.
\begin{enumerate}
\renewcommand{\labelenumi}{ \upshape{(\alph{enumi})} }
\item $1/(c+W(x)),\,c\geq0$                                                        \label{en:cone1}
\item $W^\alpha(1/x),0\leq\alpha\leq1$                                             \label{en:cone2}
\item $x^\beta W^\beta(1/x),-1\leq\beta\leq 0$                                     \label{en:cone3}
\item $W(x)/[x(c+W(x))],\,c\geq0$                                                  \label{en:cone4}
\item $1/W(x)-1/x$                                                                 \label{en:cone5}
\item $c+W(x^\beta),\,c\geq0,\:-1\leq\beta\leq 0$                                  \label{en:cone6}
\item $1/(c+W(x^\alpha)),\,c\geq0,\:0\leq\alpha\leq 1$                             \label{en:cone7}
\item $x^{\alpha\beta\gamma}W^{-\alpha\gamma}(x^\beta)[1+W(x^\beta)]^{1-\gamma},
\: 0\leq\alpha\leq1,\:-1\leq\beta\leq 0,\,0\leq\gamma\leq1$                        \label{en:cone8}
\item $1/(c+x^\alpha),\, 0\leq\alpha\leq1,\, c\geq0$                               \label{en:cone9}
\item $1-x^\alpha W^\alpha(1/x),\:0\leq\alpha\leq 1$                               \label{en:cone10}
\item $1-x^{-\alpha\beta}W^{\alpha}(x^\beta)[1+W(x^\beta)]^{-\alpha},\: 0
\leq\alpha\leq 1, \, -1\leq\beta\leq 0$                                            \label{en:cone11}
\end{enumerate}
\end{Th}
\begin{Proof}
We use the properties listed in Theorem \ref{th:StieltjesProp}.
\begin{enumerate}
\renewcommand{\labelenumi}{ \upshape{(\alph{enumi})} }
\item We apply property \rmref{en:2} to $W(x)/x$ to find that $1/W(x)\in\coneS$ and
then apply \rmref{en:3} to $1/W(x)$.
\item We first apply \rmref{en:1} to $f(x)=1/W(x)$ that is in $\coneS$ by
statement (a) and find $W(1/x)\in\coneS$.
Then we apply \rmref{en:7} to $W(1/x)$.
\item Apply \rmref{en:1} to $W(x)/x$ and apply then \rmref{en:7} to  the result.
\item Apply \rmref{en:11} to the function in the statement (a) using $W(0)=0$.
\item Apply \rmref{en:8} to $W(x)/x$ using \eqref{WnearZero} or apply \rmref{en:9} to the
function in the statement (a) with $c=0$.
\item Apply \rmref{en:5} to the function in the statement (a) and
$g(x)=x^\beta \: (-1\leq\beta\leq 0)$ that
is in $\coneS$ \cite{BergForst, Berg2008}.
\item Apply \rmref{en:4} to the function in the statement (a) and $g(x)=x^{-\alpha}\in\coneS$ for $0\leq\alpha\leq 1$.
\item Apply \rmref{en:5} to functions $f(x)=W(x)/x$ and $g(x)=x^\beta \: (-1\leq\beta\leq0)$ and find $x^\beta W^{-1}(x^\beta)\in\coneS$.
Hence by \rmref{en:7}
$a(x)=x^{\alpha\beta} W^{-\alpha}(x^\beta)\in\coneS$ for $0\leq\alpha\leq1$.
Then apply \rmref{en:5} to the function in the statement (a) with $c=1$ and $g(x)=x^\beta$ to get $b(x)=1+W(x^\beta)\in\coneS$.
Finally apply \rmref{en:6} to $a(x)$ and $b(x)$.
\item It suffices to take \rmref{en:3} with $c=1$ and $f(x)=x^{-\alpha}\, \: (0\leq\alpha\leq1)$.
\item Apply \rmref{en:12} to the function in the statement (c) with $\beta=-1$ using \eqref{WnearZero} (or apply (x) to $W(x)/x$).
\item Apply \rmref{en:10} (or \rmref{en:12}) to the result of application of \rmref{en:4} (respectively \rmref{en:5}) to the
function in the statement (d) with $c=1$ and $g(x)=x^\beta \: (-1\leq\beta\leq 0)$.
\end{enumerate}
\end{Proof}

The next theorem proves and generalizes a conjecture in \cite{Sokal}.
\begin{Th} \label{SokalFun}
The following functions are Stieltjes functions for fixed real $a\in(0,e]$
\begin{equation} \label{eq:F0}
F_0(z)=\frac{z}{1+z}W(a(1+z))/\left[W(a(1+z))-W(a)\right]^2 \ ,
\end{equation}
\begin{equation} \label{eq:F1}
F_1(z)=zW\left(\frac{a}{1+z}\right)\bigg/\left[W(a)-W\left(\frac{a}{1+z}\right)\right]^2 \ .
\end{equation}
\end{Th}
\begin{Proof}
We start with the function $F_0(z)$.
To apply Theorem \ref{th:CriterionB} we note that $F_0(z)\geq0$ for real $z>0$ ($a\in(0,e]$)
and $F_0(z)$ is a holomorphic function in the cut plane $\C\backslash(-\infty,-1-1/(ae)]$
(cf. the branch cut $\B$).
For convenience, we define a function $V(z)=\Im F_0(z)$, then it remains to show
that $V(z)\leq0$ in the upper half-plane.
Since $V(z)$ is a harmonic function in the domain $\Im z>0$, it is subharmonic there.
Thus we can apply either the maximum principle for harmonic functions in the
form of \cite[Corollary 1.10]{Axler}
or the maximum principle for subharmonic functions \cite[p.\,19-20]{Doob}.
In both cases, to get the desired result it is sufficient to ascertain that
the superior limit of $V(z)$ at all boundary points including infinity
is less than or equal to 0 \cite{AlzerBerg1}.
In other words, $V(z)\leq0$ for $\Im z>0$ if (cf. \cite[p.\,27]{Koosis})
\begin{equation*} \label{eq:cond}
\lim_{\left|z\right|\rightarrow\infty}V(z)\leq 0 \quad (\Im z>0)
\end{equation*}
and
\begin{equation} \label{eq:limsup}
\limsup_{y\rightarrow 0+}V(x+iy)\leq0 \mbox{ for all } x\in\R \ .
\end{equation}
Since $F_0(z)\sim1/\ln z$ for large $z$ owing to \eqref{eq:branch},
$V(z)\rightarrow0$ as $\left|z\right|\rightarrow\infty$
and the first condition is satisfied.
To verify the second condition we note that $V(z)$ is continuous (from above) on the
real line $z=x\in\R$.
Then after introducing variable $t=a(1+x)$ the inequality \eqref{eq:limsup} becomes
\begin{equation*} \label{ineq:Im}
\Im \left\{\left(1-\frac{a}{t}\right)\frac{W(t)}{\left[W(t)-W(a)\right]^2}\right\}\leq0
\end{equation*}
\\
or substituting $W(t)=u+iv$ and using \eqref{eq:u(v)}

\begin{equation} \label{ineq:v}
\frac{-v}{\left[(b+v\cot v)^2+v^2\right]^2}\left(\frac{v^2}{\sin^2v}-b^2\right)
\left(1-\frac{a}{t}\right)\leq0 \ ,
\end{equation}
\\
where $b=W(a)\in(0,1]$.
\\We have $v=0$ for $t\geq-1/e$ and $0<v<\pi$ for $-\infty<t<-1/e$ (see \eqref{eq:ImW}).
Since $0<b\leq1$ and $v^2/\sin^2v>1$ for $v\in(0,\pi)$, the inequality \eqref{ineq:v} holds
for all real $t$. Thus $F_0(z)$ is a Stieltjes function.

The theorem for the function $F_1(z)$ follows from the relation
\begin{equation} \label{rel:F0F1}
F_1(z)=-F_0\left(-\frac{z}{1+z}\right)
\end{equation}
because in terms of the conditions of Theorem \ref{th:CriterionB} the transformation
in the right hand-side of \eqref{rel:F0F1} retains the properties of $F_0(z)$.
In particular, $\Im F_1(z)\leq0$ for $\Im z>0$ because, first,
 $\Im z$ and $\Im(-z/(1+z))$ are of the opposite signs
and secondly, $\Im F_0(z)\geq0$ for $\Im z<0$
which follows from $F_0(\bar{z})=\overline{F_0(z)}$ due to near conjugate symmetry and
the established above non-positivity of $\Im F_0(z)$ in the upper half-plane.
Thus $F_1(z)$ is also a Stieltjes function.
\end{Proof}
\begin{remark}
We make a note about a behavior of functions \eqref{eq:F0} and \eqref{eq:F1} for large and small $z$.
Specifically, using \eqref{eq:branch} and \eqref{WnearZero} one can obtain respectively
$F_0(z)\rightarrow0$ and $F_1(z)\rightarrow a/W^2(a)$ as $z\rightarrow\infty$.
Using \eqref{eq:DW} we find $F_{0,1}\sim c/z$ as $z\rightarrow0$, where $c=(1+W(a))^2/W(a)$.
By the way, owing to property \rmref{en:9} given in theorem \ref{th:StieltjesProp}
the differences $F_{0,1}-c/z$ are Stieltjes functions as well.
\end{remark}

\subsection{Is $W$ a Stieltjes function?} \label{sec:W-S}
The principal branch of the Lambert $W$ function itself is not a Stieltjes function
in the sense of Definition \ref{def:Stieltjes fun}. This can be shown in different ways. For example,
one can apply~Theorem \ref{th:criterionA} to $W(z)$ to see that the second condition \eqref{cond:f} fails.
Indeed, when $z=is$ we have by \eqref{Imz:u(v)} and \eqref{Imz:s(v)}
\begin{equation*}
\left|sW(is)\right|=s\sqrt{u^2+v^2}=v^2\sec^2(v)e^{v\tan v}\rightarrow\infty \quad \mbox{as} \quad v\rightarrow\pi/2.
\end{equation*}
The same conclusion can be reached using Theorem \ref{th:CriterionB} because \eqref{signImW}
contradicts \eqref{inequalities}.
Finally, $W$ is not a Stieltjes function because it is not an anti-Herglotz function
(cf. Remark \ref{rem:Herglotz}).

Note, however, that $W$ function can be regarded as a Stieltjes function
in the sense of a definition given in \cite{Tokarzewski} and \cite{Brodsky} or
used in \cite{Telega} and different from \eqref{eq:DefStieltFun} by the factor $z$ in
the right hand-side.
$W$ function can also be considered as a generalized Stieltjes transform by the
definition in \cite{SaxenaGupta} (which is different from that of the generalized
Stieltjes transform defined in \cite[p.\,30]{WidderSt} and studied, for example,
in \cite{Schwarz} and \cite{Sokal:Stieltjes}).
Finally,
in~\cite{Schilling}, the terms Stieltjes function and Stieltjes representation are not treated as
equivalent (compare definitions~\cite[p.\,11]{Schilling} and \cite[p.\,55]{Schilling}).
By these definitions $W(z)$ has a Stieltjes representation
(which is the result of multiplication of the representation \eqref{eq:result2} by $z$)
though it is not a Stieltjes function.

\section{Particular Stieltjes integrals}
The Stieltjes representation for $W/z$ given in \eqref{eq:result2} and \eqref{eq:measure}
itself contains $W$, which can be regarded as self-referential. Here we give representations
containing only elementary functions.

\begin{Th} \label{th:W/z}
The following representation of function $W(z)/z$ holds \cite{poster}

\begin{equation} \label{IntegralW/z}
\frac{W(z)}{z}=\frac{1}{\pi}\int_0^\pi\frac{v^2+(1-v\cot v)^2}{z+v\csc(v)e^{-v\cot v}}dv \quad (\left|\arg z\right|<\pi) \ .
\end{equation}
\end{Th}
\begin{Proof}
We start with \eqref{eq:res1} and, noting \eqref{eq:ImW}, change to the variable
$v=\Im W(t)$. The integral becomes
\begin{equation} \label{W(z)/z:Int(v)}
\frac{W(z)}{z}=\frac{1}{\pi}\int_0^\pi\frac{v}{t(z-t)}\frac{dv}{v^\prime(t)} \ ,
\end{equation}
where the variables $t$ and $v$ are related by \eqref{eq:t(v)}.
The derivative $v^\prime(t)$ is conveniently found by taking the imaginary part
of \eqref{eq:DW}, and using \eqref{eq:u(v)}.
\begin{equation} \label{eq:v'(t)}
v^\prime(t)=\frac{v}{t\left[v^2+(1+u)^2\right]}=\frac{v}{t\left[v^2+(1-v\cot v)^2\right]} \ .
\end{equation}
\end {Proof}

\begin{remark}
Since the integrand in \eqref{IntegralW/z} is an even function (with respect to $v$),
the integral admits the symmetric form
\begin{equation}
\frac{W(z)}{z}=\frac{1}{2\pi}\int_{-\pi}^\pi\frac{v^2+(1-v\cot v)^2}{z+v\csc(v)e^{-v\cot v}}dv \quad (\left|\arg z\right|<\pi) \ . \notag
\end{equation}
This integral has a $C^\infty$ periodic extension and thus the midpoint rule
is spectrally convergent for its quadrature (see e.g. \cite{Weideman}).
\end{remark}

We now pay attention to the statement (d) of Theorem \ref{th:coneS} with $c=1$ which means
by \eqref{eq:DW} that the derivative $W^\prime(x)\in\mathcal{S}$.
We derive an integral representation of $W^\prime(z)$ in the complex $z$-plane in the next theorem.

\begin{Th} \label{th:derivativeW}
The derivative of $W$ function has the following Stieltjes-integral representation
\begin{equation} \label{IntegralW'}
W^\prime(z)=\frac{1}{\pi}\int_0^\pi\frac{dv}{z+v\csc(v)e^{-v\cot v}} \quad (\left|\arg z\right|<\pi) \ .
\end{equation}
\end{Th}
\begin{Proof}
We take the formula \eqref{def:complexSt} with $a=0$ due to \eqref{DWnearInfty}.
After changing $t$ by $-t$ and introducing the function
$\rho(t)=-\mu(-t)$ we obtain
\begin{equation} \label{Int1}
W^\prime(z)=\int_{-\infty}^0\frac{d\rho(t)}{z-t} \ .
\end{equation}
The function $\rho(t)$ can be determined using the Stieltjes-Perron inversion formula \cite[p.\,591]{Henrici2}
\begin{equation*}
\rho(t)=-\frac{1}{\pi}\lim_{s\rightarrow 0+}\Im\int_{-\infty}^tW^\prime(t+is)dt \
\end{equation*}
for all continuity points $t$.
Since $\rho(t)$ is defined to arbitrary constant, after integrating one can set
\begin{equation} \label{def:rho(t)}
\rho(t)=-\frac{1}{\pi}\lim_{s\rightarrow 0+}\Im W(t+is)= -\frac{1}{\pi}\Im W(t)\ ,
\end{equation}
where the limit uses the continuity from above of $W$ on its branch cut.
The same result can be obtained using one of Sokhotskyi's formulas \cite[p.\,138]{Henrici3}.

Since $\rho(t)$ increases from $-1$ to $0$ as $t$ changes
from $-\infty$ to $-1/e$ and vanishes on the
remainder of the real line (cf. \eqref{eq:ImW}), the domain of integration in \eqref{Int1} is
defined by $-\infty<t<-1/e$. Owing to the relation
\begin{equation} \label{mu(rho)}
\mu(t)=-\rho(-t)\ ,
\end{equation}
the function $\mu(t)$ in \eqref{def:realSt} can be regarded as a positive measure.
Thus \eqref{Int1} takes the form
\begin{equation} \label{Int2}
W^\prime(z)=\frac{1}{\pi}\int_{-\infty}^{-1/e}\frac{1}{t-z}\frac{d\Im W(t)}{dt}dt \ .
\end{equation}\\
This formula can also be found by considerations similar to those used in
the proof of Theorem \ref{th:W(z)/z}.
By the change of variables $v=\Im W(t)$ in the integral \eqref{Int2} (see
the proof of Theorem \ref{th:W/z} and relation \eqref{eq:t(v)}) we obtain \eqref{IntegralW'}.
\end{Proof}
\begin{corollary}
\begin{equation} \label{eq:Nutend}
\int_0^\pi\left\{\frac{\sin v}{v}e^{v\cot v}\right\}^\nu dv=\frac{\pi \nu^\nu}{\nu!},  \quad  \nu\in\N \ .
\end{equation}
\end{corollary}
\begin{Proof}
The integral \eqref{IntegralW'} can be written as
\begin{equation} \label{eq:Nutstart}
\sum_{n=1}^\infty(-1)^{n-1}\frac{n^n}{n!}z^{n-1}=\frac{1}{\pi}\int_0^\pi\frac{dv}{z-t} \ ,
\end{equation}\\
where $t$ is defined by \eqref{eq:t(v)} and the left-hand side is obtained by differentiation of
the series \eqref{eq:Taylor} that is convergent for $\left|z\right|<1/e$. Since $\left|t\right|>1/e$ and
therefore $\left|z\right|<\left|t\right|$, we can expand $(z-t)^{-1}$ in the non-negative powers of $z$.
Equating the coefficients of the same power of $z$ in \eqref{eq:Nutstart} we obtain an equality
\[
(-1)^{n-1}\frac{n^n}{n!}=-\frac{1}{\pi}\int_0^\pi\frac{dv}{t^n}
\]
which after substituting \eqref{eq:t(v)} results in \eqref{eq:Nutend}.
\end{Proof}

It is obvious that if the integral  \eqref{eq:Nutend} is known then going back from it
to \eqref{eq:Nutstart} we find \eqref{IntegralW'}.
The integral \eqref{eq:Nutend} was conjectured by Nuttall for $\nu\geq0$ \cite{Nuttall};
Bouwkamp found a more general
integral \cite{Bouwkamp}, for which Nuttall's conjecture is a special case,
using a representation of $\pi \nu^\nu/\Gamma(\nu+1)$ via a Hankel-type integral.
Thus a representation of the derivative of $W$ function in the form of the
Stieltjes integral \eqref{IntegralW'} allows one to compute the integral \eqref{eq:Nutend} and conversely,
starting with the integral of Nuttall--Bouwkamp one can obtain formula \eqref{IntegralW'} in a way
completely different from that used in the proof of Theorem \ref{th:derivativeW}.
It is interesting to note that the connection between \eqref{eq:Nutend} and Lambert $W$ was
noted by W.E. Hornor and C.C. Rousseau before $W$ was named (see editorial remarks in \cite{Nuttall}).

Coming back to the results of Theorem \ref{th:coneS} we consider the assertion
(a) with $c=1$ and assertion (e) by which $1/(1+W(z))\in\coneS$
and $1/W(z)-1/z\in\coneS$.
We can derive integral representations of these functions in the same manner
as it was done for $W^\prime(z)$ in the proof of Theorem \ref{th:derivativeW}.
The result is in the following theorem.

\begin{Th}
The following Stieltjes-integral representations hold
\begin{equation} \label{eq:1/1+W}
\frac{1}{1+W(z)}=\frac{1}{\pi}\int_0^\pi\frac{dv}{1+ze^{v\cot v}\sin v/v}
\quad (\left|\arg z\right|<\pi) \ ,
\end{equation}

\begin{equation} \label{eq:1/W}
\frac{1}{W(z)}=\frac{1}{z}+\frac{1}{\pi}\int_0^\pi\frac{v^2+(1-v\cot v)^2}{v\csc(v)\left(v\csc(v)+ze^{v\cot v}\right)}dv
\quad (\left|\arg z\right|<\pi) \ .
\end{equation}
\end{Th}
\vspace{5mm}
\begin{corollary}
\begin{equation} \label{eq:W=ln}
W(z)=\ln\left[1+\frac{z}{\pi}\int_0^\pi\frac{v^2+(1-v\cot v)^2}{v\csc(v)\left(v\csc(v)+ze^{v\cot v}\right)}dv\right] \ .
\end{equation}
\end{corollary}
\begin{Proof}
By substituting \eqref{eq:1/W} in $W(z)=\ln(z/W(z))$ (cf. \eqref{eq:Wright}).
\end{Proof}
\section{Completely monotonic functions} \label{sec:cmf}
We denote by $\mathcal{CM}$ the set of all completely monotonic functions, which are
defined as follows~\cite{AlzerBerg}.
\begin{definition} A function $f:(0,\infty)\rightarrow\R$ is called a \textit{completely monotonic} function
if $f$ has derivatives of all orders and satisfies
$(-1)^nf^{(n)}(x)\geq0$ for $x>0$, $n=0,1,2,...$
\end{definition}
The set of Stieltjes functions is contained in the set of completely monotonic functions, and thus
all of the functions listed in Theorem \ref{th:coneS} are completely monotone.
The set $\mathcal{CM}$ is a convex cone containing the positive constant functions;
a product of completely monotonic functions is again completely monotone \cite[p.\,61]{BergForst}.
By Bernstein's theorem
\cite[Theorem 9.3]{BergForst}, a function $f\in\mathcal{CM}$ if and only if it is of the form
\begin{equation} \label{eq:monotone}
f(x)=\int_0^\infty e^{-x\xi}d\nu(\xi) \quad (x>0),
\end{equation}
where $\nu$ is an uniquely determined positive measure on $[0,\infty)$.
Completely monotonic functions are in turn connected with the set of Bernstein functions.

\begin{definition} \cite[Definition 5.1]{Berg2008} A function $f: (0,\infty)\rightarrow[0,\infty)$
is called a \textit{Bernstein function} if it is $C^\infty$ and $f^\prime$ is completely monotonic.
\end{definition}
We denote the set of Bernstein functions by $\mathcal{B}$.
Since $W^\prime\in\coneS\subset\mathcal{CM}$, $W$ is a Bernstein function.
A Bernstein function $f(x)$ admits the L\'evy-Khintchine representation
\begin{equation} \label{def:Bernstein}
f(x)=a+bx+\int_0^\infty\left(1-e^{-x\xi}\right)d\nu(\xi) \ ,
\end{equation}
where $a,b\geq0$ and $\nu$ is a positive measure on $(0,\infty)$ satisfying
$\int_0^\infty \xi(1+\xi)^{-1}d\nu(\xi)<\infty$. It is called the L\'evy measure.
The equation \eqref{def:Bernstein} is obtained by integrating \eqref{eq:monotone} written for $f^\prime$ \cite{Berg2008}.

An important relation between the classes $\coneS$ and $\mathcal{B}$ is given by the assertion \cite[Theorem 5.4]{Berg2008}
\begin{equation} \label{rel:SB}
g\in\coneS\:\backslash\left\{0\right\}\Rightarrow 1/g\in\mathcal{B}.
\end{equation}
Combining this with the function composition result \cite[Corollary 5.3]{Berg2008} that
$f\in\mathcal{CM}$ and $g\in\mathcal{B}$ implies $f\circ g\in\mathcal{B}$, we obtain
the following lemma
\begin{Lem} \label{lemmma:CM}
If $f\in\mathcal{CM}$ and $g\in\coneS\:\backslash\left\{0\right\}$ then $f(1/g)\in\mathcal{CM}$.
\end{Lem}
This lemma extends the list of completely monotonic functions containing $W$.
\begin{Th}
The following functions are completely monotonic
\begin{enumerate}
\renewcommand{\labelenumi}{ \upshape{(\alph{enumi})} }
\item $x^\lambda W(x)$ ($x>0, \: \lambda\leq-1$).
\item $x^\lambda W^\alpha(x^\beta)\left[1+W(x^\beta)\right]^\gamma \quad (x>0, \:  \alpha,\gamma\geq0, \, -1\leq\beta\leq 0, \, \lambda\leq0)$.
\item $x^\lambda W^\alpha(x^{-\beta})\left[1+W(x^{-\beta})\right]^\gamma \quad (x>0, \: \alpha,\gamma\leq0, \, -1\leq\beta\leq 0, \, \lambda\leq0)$.
\item $1-x^{-\alpha\beta\gamma}W^{\alpha\gamma}(x^\beta)[1+W(x^\beta)]^{\gamma-1} \quad (x>0, \:  0\leq\alpha\leq 1, \,
-1\leq\beta\leq 0, \, 0\leq\gamma\leq 1)$.
\end{enumerate}
\end{Th}
\begin{Proof}
\begin{enumerate}
\renewcommand{\labelenumi}{ \upshape{(\alph{enumi})} }
\item Since $W(x)/x\in\coneS\subset\mathcal{CM}$ and $x^\alpha\in\mathcal{CM}$ for $\alpha\leq0$, the function $x^\lambda W(x)$ ($\lambda\leq-1$) is a product of two completely monotonic functions and the statement (a) follows.
\item Take function $f_\alpha(x)=x^{-\alpha}\in\mathcal{CM}$ ($x>0, \, \alpha\geq0$) and functions $g(x)=1/W(x^\beta)$ and $h(x)=1/(1+W(x^\beta))$ where $-1\leq\beta\leq0$. Since $1/g\in\coneS$ and $1/h\in\coneS$ by Theorem \ref{th:coneS} (f) with $c=0$ and $c=1$ respectively, by Lemma \ref{lemmma:CM} we have $f_\alpha(g(x))=g^{-\alpha}(x)\in\mathcal{CM}$ and $f_\gamma(h(x))=h^{-\gamma}(x)\in\mathcal{CM}$ ($\gamma\geq0$).
Substituting functions $g(x)$ and $h(x)$ in the power functions and taking a product of obtained completely monotonic functions with $x^\lambda\in\mathcal{CM}$ ($x>0, \, \lambda\leq0$), the statement (b) follows.
\item Consider function $f_\lambda(x)=x^\lambda \in\mathcal{CM}$ ($x>0, \, \lambda\leq0$) and functions $g(x)=W(x^{-\beta})$ and $h(x)=1+W(x^{-\beta})$ where $-1\leq\beta\leq0$. Since $1/g\in\coneS$ and $1/h\in\coneS$ by Theorem \ref{th:coneS} (g) with $c=0$ and $c=1$ respectively, by Lemma \ref{lemmma:CM} we have $f_\alpha(g(x))=g^\alpha(x)\in\mathcal{CM}$ and $f_\gamma(h(x))==h^\gamma(x)\in\mathcal{CM}$ for $\alpha\leq0$ and $\gamma\leq0$.
Substituting functions $g(x)$ and $h(x)$ and taking a product of obtained functions with $f_\lambda(x)$, the statement (c) follows.
\item By Theorem \ref{th:coneS} (h) and  the assertion \eqref{rel:SB},
for $x>0, \, 0\leq\alpha\leq 1, \, -1\leq\beta\leq 0, \, 0\leq\gamma\leq 1$ we have
$f(x)=g^{\alpha\gamma}(x)[1+W(x^\beta)]^{\gamma-1}\in\mathcal{B}$,
where $g(x)=x^{-\beta}W(x^\beta)$. In addition, the function $f(x)$ is bounded,
particularly, $0<f(x)<1$ because $0<[1+W(x^\beta)]^{\gamma-1}<1$ and $0<g(x)<1$
(the latter follows from the fact that $g(x)$ goes to 0 and 1 as $x$ tends to 0 and $\infty$ respectively and $g^\prime(x)>0$,
which can be established using \eqref{WnearZero}, \eqref{WnearInfty} and \eqref{eq:DW}).
Then by \cite[Remark 5.5]{Berg2008} the assertion (d) follows.
\end{enumerate}
\end{Proof}

We considered only sufficient conditions for a function to be a completely monotonic.
To find the necessary and sufficient conditions is a much more complicated problem so that
in some cases it requires (at least as the first step) using
the methods of experimental \mbox{mathematics \cite{JeffreyMono}}.

\section{Complete Bernstein functions} \label{sec:cbf}
A very important subclass in $\mathcal{B}$ is the class of complete Bernstein functions
 denoted by $\mathcal{CB}$.
\begin{definition}\cite[Definition 6.1]{Schilling} A Bernstein function $f$ is called
a \textit{complete Bernstein function} if the L\'evy measure
in \eqref{def:Bernstein} is such that $d\nu(t)/dt$ is a completely monotonic function.
\end{definition}
We point out four connections between classes $\mathcal{CB}$ and $\coneS$ used in this paper
(for additional relations between these classes see \cite[Chapter 7]{Schilling}).
By Proposition 7.7 in \cite{Schilling},
\begin{equation} \label{connect1:CB,S}
f\in\coneS\Rightarrow f(0)-f(x)\in\mathcal{CB} \ ,
\end{equation}
where the limit of $f(x)$ at $x=0$ (from the right) is assumed to be finite.
Also if $f$ is bounded and $f\in \mathcal{CB}$, there exists a bounded $g\in\coneS$ with
$\lim_{x\rightarrow\infty}g(x)=0$ such that
\begin{equation} \label{connect2:CB,S}
 f(x)=f(0)+g(0)-g(x) \ .
\end{equation}
In addition, \cite[Theorem 7.3]{Schilling} and \cite[Theorem 6.2(i),(ii)]{Schilling} establish
\begin{align} \label{connect3:CB,S}
f&\in\mathcal{CB}\Leftrightarrow 1/f\in\coneS\:\backslash\left\{0\right\} \ ,\\
 \label{connect4:CB,S}
f&\in\mathcal{CB}\Leftrightarrow f(x)/x\in\coneS \ .
\end{align}

Now we go back to the properties of the set $\coneS$ listed in
section \ref{sec:properties} to prove the last three properties therein.
Let $f\in\coneS\:\backslash\left\{0\right\}$.

(x) Apply sequentially (vii), \eqref{connect1:CB,S}, \eqref{connect3:CB,S}, (i), to obtain
$f^\alpha\in\coneS\;(0\leq\alpha\leq1)\Rightarrow f^\alpha(0)-f^\alpha(x)\in\mathcal{CB}
\Rightarrow g(x)=\left[f^\alpha(0)-f^\alpha(x)\right]^{-1}\in\coneS\Rightarrow 1/g(1/x)=f^\alpha(0)-f^\alpha(1/x)\in\coneS$;

(xi) Apply sequentially \eqref{connect1:CB,S}, \eqref{connect3:CB,S}, (ii), to obtain
$f(0)-f(x)\in\mathcal{CB}\Rightarrow
g(x)=\left[f(0)-f(x)\right]^{-1}\in\coneS\Rightarrow 1/(xg(x))=(f(0)-f(x))/x\in\coneS\Rightarrow(1-f(x)/f(0))/x\in\coneS$;

(xii) By (vii), $f^\alpha\in\coneS\;(0\leq\alpha\leq1)$. Suppose that $\lim_{x\rightarrow0}f(x)=b\leq\infty$
and $\lim_{x\rightarrow\infty}f(x)=c$ where $0<c<\infty$.
Then $b^{-\alpha} \leq f^{-\alpha} \leq c^{-\alpha}$, i.e. $f^{-\alpha}$ is bounded.
In addition, $f^{-\alpha}\in\mathcal{CB}$ by \eqref{connect3:CB,S}.
Therefore the statement \eqref{connect2:CB,S} can be applied, i.e. there exists a bounded function
$g\in\coneS, \lim_{x\rightarrow\infty}g(x)=0$ such that we can write $g(x)=g(0)+b^{-\alpha}-f^{-\alpha}(x)$.
Taking the last equation in the limit $x\rightarrow\infty$ we obtain $g(0)+b^{-\alpha}=c^{-\alpha}$,
hence $g=c^{-\alpha}-f^{-\alpha}$ and the assertion follows.

In closing this section we note that the statement \eqref{connect3:CB,S}
with $1/W \in\coneS$ (by Theorem \ref{th:coneS}(a) with $c=0$) immediately results in $W\in\mathcal{CB}$.
Being a complete Bernstein fucntion $W$ has an integral representation
that is the result of multiplication
of \eqref{IntegralW/z} by $z$ \cite[Remark 6.4]{Schilling},
which reflects the relation \eqref{connect4:CB,S}.
In addition, the complete Bernstein functions are closely connected to the Pick functions considered in Section \ref{sec:Pick}.

\section{Bernstein representations}
Not only does $W\in\mathcal{CB}$ as shown, it also belongs to another subset of Bernstein functions.

\begin{definition}\label{def:TBF} \cite[Definition 8.1]{Schilling} A Bernstein function $f$
is called a \textit{Thorin-Bernstein function} if the L\'evy measure in \eqref{def:Bernstein} is
such that $t\,d\nu(t)/dt$ is a completely monotonic function.
\end{definition}
To find out whether $W$ is a Thorin-Bernstein function
we apply Theorem 8.2 in \cite{Schilling}, which establishes five equivalent assertions (i)-(v)
which we refer to below.
In particular, in accordance with  assertions (i) and (ii),
$W(x)$ is a Thorin-Bernstein function because $W(x)$ maps $(0,\infty)$ to itself,
$W(0)=0$ and $W^\prime(x)\in\coneS$.
As a Thorin-Bernstein function, $W(x)$ admits the following integral representation.

\begin{Th}
The principal branch of the $W$ function can be represented as the integral
\begin{equation} \label{eq:Thorin}
W(z)= \frac{1}{\pi}\int_0^\pi\ln\left(1+z\frac{\sin v}{v}e^{v\cot v}\right)dv \quad (\left|\arg z\right|<\pi) \ .
\end{equation}
\end{Th}

\begin{Proof}
It follows from the proof of Theorem 8.2 in \cite{Schilling} and its assertion (v) that
$W(x)$ being a Thorin-Bernstein function, it can also be written in the form
\begin{equation} \label{eq:Thorin2}
W(x)=a+bx+\int_0^\infty \frac{x}{x+t}\frac{\tau(t)}{t}dt \ ,
\end{equation}
where due to $W(0)=0$ and \eqref{WnearInfty} we have $a=0$ and $b=0$.
Comparing \eqref{eq:Thorin2} and \eqref{eq:result2}, \eqref{eq:measure} we find $\tau(t)=\Im W(-t)/\pi$.
Integration by parts gives
\begin{equation} \label{eq:Thorin3}
W(x)=-\frac{1}{\pi}\int_{-\infty}^{-1/e} \ln\left(1-\frac{x}{t}\right)\frac{d}{dt}\Im W(t)dt \ ,
\end{equation}
where we changed $t\rightarrow -t$.
Changing the variables $v=\Im W(t)$ in \eqref{eq:Thorin3} (see \eqref{eq:t(v)}) and taking
a holomorphic  extension of the result to the cut $z$-plane $\C\backslash(-\infty,0]$
satisfying near conjugate symmetry, we obtain \eqref{eq:Thorin}.
\end{Proof}

\begin{remark}
Differentiating the representation \eqref{eq:Thorin} for $W(z)$ gives formula \eqref{IntegralW'} for $W^\prime(z)$.
\end{remark}
\begin{remark}
The representation \eqref{eq:Thorin3} was obtained in \cite{Caillol} as a dispersion relation for
the principal branch of $W$ function using Cauchy's integral formula in a manner similar to the method applied for the proof of Theorem \ref{th:W(z)/z}.
\end{remark}
As a Bernstein function, $W$ can be written in the form \eqref{def:Bernstein}, with $a=b=0$ as in
\eqref{eq:Thorin2},
allowing us to establish one more representation of $W$.

\begin{Th}
For the principal branch of the $W$ function the following formula holds
\begin{equation} \label{eq:Bernstein}
W(z)= \int_0^\infty\frac{1-e^{-z\xi}}{\xi}\varphi(\xi)d\xi \quad (\Re z\geq0),
\end{equation}
where
\begin{equation} \label{eq:varphi}
\varphi(\xi)=\frac{1}{\pi}\int_0^\pi\exp\left(-\xi v\csc(v)e^{-v\cot v}\right)dv.
\end{equation}
\end{Th}
\begin{Proof}
We consider the Stieltjes integral \eqref{def:realSt} for the derivative
\begin{equation*}
W^\prime(x)=\int_0^\infty\frac{d\mu(\theta)}{x+\theta}
\end{equation*}
and use the representation
$(x+\theta)^{-1}=\int_0^\infty e^{-(x+\theta)\xi}d\xi$
to write it in the form
\begin{equation} \label{DoubleInt}
W^\prime(x)=\int_0^\infty\left\{\int_0^\infty e^{-\xi\theta}d\mu(\theta)\right\}e^{-x\xi}d\xi \ .
\end{equation}
Comparing \eqref{DoubleInt} and the result of differentiating \eqref{def:Bernstein} we find the relation between measures $\mu$ and $\nu$ \cite{BergGamma}
\begin{equation*} \label{eq:measures}
\frac{d\nu}{d\xi}=\frac{1}{\xi}\int_0^\infty e^{-\xi\theta}d\mu(\theta)\ .
\end{equation*}
Using the relation \eqref{mu(rho)} and formula \eqref{def:rho(t)},
after changing the variables $v=\Im W(-\theta)$ (see \eqref{eq:t(v)}) we obtain
\begin{equation} \label{dnu/dt}
d\nu=\frac{\varphi(\xi)}{\xi}d\xi     \ ,
\end{equation}
where $\varphi(\xi)$ is defined by \eqref{eq:varphi}.
We collect the intermediate results and take a holomorphic continuation of \eqref{def:Bernstein}
to the right half-plane $\Re z\geq0$,
where the integral  \eqref{eq:Bernstein} is convergent, in accordance with
near conjugate symmetry (cf. Proposition 3.5 in \cite{Schilling}).
\end{Proof}
Note that by \eqref{eq:varphi} function $\varphi(\xi)\in\mathcal{CM}$, as should be, because $W$ is still a Thorin-Bernstein function
(cf. Definition \ref{def:TBF})

\section{Pick representations} \label{sec:Pick}
\begin{definition} \cite[Definition 4.1]{Berg2008}
A function $f(z)$ is called a \textit{Pick function} (or \textit{Nevanlinna function}) if it is holomorphic
 in the upper half-plane $\Im z>0$ and $\Im f\geq0$ there.
\end{definition}
A Pick function $f(z)$ admits an integral representation \cite[Theorem 4.4]{Berg2008}
\begin{equation} \label{int:genPick}
f(z)=\alpha_0+b_0z+\int_{-\infty}^\infty\frac{1+tz}{(t-z)(1+t^2)}d\sigma(t) \quad (\Im z > 0) \ ,
\end{equation}
where
\begin{equation} \label{Pick:coeff}
\alpha_0=\Re f(i), \quad b_0=\lim_{y\rightarrow\infty}\frac{f(iy)}{iy} \ ,
\end{equation}
and a positive measure $\sigma$ satisfies
\begin{equation} \label{Pick:measure}
\lim_{s\rightarrow0+}\frac{1}{\pi}\int_\R\Im f(t+is)\varphi(t)dt=\int_\R\varphi(t)d\sigma(t) \
\end{equation}
for all continuous functions $\varphi:\R\rightarrow\R$ with compact support.
The formula \eqref{int:genPick} with the integral written in terms of a measure $d\tilde{\sigma}(t)=\pi(1+t^2)^{-1}d\sigma(t)$
is called a \textit{Nevanlinna formula} \cite[p.\,100]{Levin}.

Since $W(z)$ is a holomorphic function in the upper half-plane $\Im z > 0$ with the property \eqref{signImW},
$W(z)$ is a Pick function. It also follows from the two facts that $W\in\mathcal{CB}$ (see Section \ref{sec:cbf})
and that the complete Bernstein functions are exactly those Pick functions
which are non-negative on the positive real line \cite[Theorem 6.7]{Schilling}.
Thus $W$ admits a representation \eqref{int:genPick} and in view of that the following theorem holds.

\begin{Th}
The principal branch of $W$ function can be represented in the form
\begin{equation} \label{Pick1}
W(z)=\alpha_0+\frac{1}{\pi}\int_0^\pi K(z,v)t(v)dv \quad (\left|\arg z\right|<\pi) \ ,
\end{equation}
where $\alpha_0=\Re W(i)=0.3746990..$,
\begin{equation} \label{K(z,v)}
K(z,v)=\frac{\left(1+zt(v)\right)\left(v^2+(1-v\cot v)^2\right)} {(z-t(v)) \left(1+t^2(v)\right)} \ ,
\end{equation}\\
and $t(v)$ is defined by \eqref{eq:t(v)}.
\end{Th}

\begin{Proof}
Apply formulas \eqref{Pick:coeff}-\eqref{Pick:measure} to function $f(z)=W(z)$
\[
\alpha_0=\Re W(i), \quad b_0=\lim_{y\rightarrow\infty}\frac{W(iy)}{iy},
\quad d\sigma(t)=\frac{1}{\pi}\Im W(t)dt.
\]
Using \eqref{WnearInfty}, we see $b_0=0$.
Since $\Im W(t)=0$ for $t\geq-1/e$ (cf. \eqref {eq:ImW}), we obtain
\begin{equation} \label{eq:PickRes}
W(z) = \alpha_0 + \frac{1}{\pi}\int_{-\infty}^{-1/e}\frac{1+tz}{(t-z)(1+t^2)}\Im W(t)dt \quad (\Im z > 0) \ .
\end{equation}
By the change of variables $v=\Im W(t)$ in the integral \eqref{eq:PickRes} (see \eqref{eq:t(v)})
we obtain formula
\eqref{Pick1} that is also valid in the lower half-plane $\Im z<0$ in accordance with
near conjugate symmetry of $W$.
\end{Proof}
\begin{corollary}
\begin{equation} \label{Pick2}
\frac{W(z)}{z}=\gamma_0\exp\left\{-\frac{1}{\pi}\int_0^\pi K(z,v)t(v)dv \right\} \quad (\left|\arg z\right|<\pi) \ ,
\end{equation}\\
where $\gamma_0=e^{-\Re W(i)}=0.6874961..$.
\end{corollary}
\begin{Proof}
It immediately follows from \eqref{Pick1} owing to the identity $W(z)/z=e^{-W(z)}$.
\end{Proof}
Now we take advantage of the fact that if function $f\in\coneS$ then $-f$ and $1/f$ are
Pick functions \cite{Berg2008}.
Therefore, since $W(x)/x\in\coneS$, $-W(x)/x$ and $x/W(x)$ are Pick functions that admit a
representation \eqref{int:genPick}.
We can obtain a representation \eqref{int:genPick} for functions $-W(x)/x$ and $x/W(x)$ similar
to the derivation of formula \eqref{Pick1}, and the result is in the following theorem.

\begin{Th}
For the principal branch of the $W$ function the following formulas hold
\begin{equation} \label{Pick3}
\frac{W(z)}{z}=\beta_0+\frac{1}{\pi}\int_0^\pi K(z,v)dv \quad (\left|\arg z\right|<\pi) \ ,
\end{equation}
\begin{equation} \label{Pick4}
\frac{z}{W(z)}=\eta_0-\frac{1}{\pi}\int_0^\pi K(z,v)e^{-2v\cot v}dv \quad (\left|\arg z\right|<\pi) \ ,
\end{equation}\\
where $K(z,v)$ is defined by \eqref{K(z,v)}, $\beta_0=\Re\left[W(i)/i\right]=\Im W(i)=0.5764127..$,
$\eta_0=\Re[i/W(i)]=1.21953\ldots$
\end{Th}

The constants in \eqref{Pick1} and \eqref{Pick2}-\eqref{Pick4} obey the
relations $\alpha_0+i\beta_0=W(i)$, $\gamma_0=e^{-\alpha_0}=\beta_0/\cos\beta_0$,
$\eta_0=\beta_0/(\alpha_0^2+\beta_0^2)$.

We add in one more integral representation associated with the Nevanlinna formula
which follows from the result obtained by Cauer \cite{Cauer}.
Specifically, based on the Riesz-Herglotz formula \cite[p.\,99]{Levin}
Cauer proved that if a real symmetric function $f(z)$ with non-negative real part is
holomorphic in the right $z$-half-plane,
it can be represented as
\begin{equation} \label{Poisson:f(z^2)}
f(z)=z\left[b+\int_0^\infty\frac{dh(r)}{z^2+r}\right] \ ,
\end{equation}
where constant $b\geq0$ and
\begin{equation} \label{def:h(r)}
h(r)=\frac{2}{\pi}\lim_{x\rightarrow0}\Re\int_0^{\sqrt{r}} f(x+iy)dy \ .
\end{equation}
In fact, the formula \eqref{Poisson:f(z^2)} follows from the Nevanlinna formula (or \eqref{int:genPick})
after changing the variable $z\rightarrow-iz$, which transforms the upper half-plane onto the right half-plane,
and taking into account $f(\bar{z})=\overline{f(z)}$.
\begin{Th}
The following representation of function $W(z)/z$ holds
\\
\begin{equation} \label{int:z^2}
\frac{W(z)}{z}=\frac{2}{\pi}\int_0^{\pi/2}\frac{\left[v^2+(1+v\tan v)^2\right]v\sec(v)e^{v\tan v} }
{z^2+v^2\sec^2(v)e^{2v\tan v}}\tan v\,dv \quad (\Re z>0)\ .
\end{equation}
\end{Th}
\begin{Proof}
Since $W$ function meets the above requirements, the formulas \eqref{Poisson:f(z^2)} and
\eqref{def:h(r)} can be applied with the result
\begin{equation*}
\frac{W(z)}{z}=\frac{2}{\pi}\int_0^\infty\frac{\Re W(is)}{z^2+s^2}ds \quad (\Re z>0)\ ,
\end{equation*}
\\where we set $b=0$ due to \eqref{WnearInfty} and $r=s^2$.
\\Changing the variables defined by \eqref{Imz:u(v)}-\eqref{Imz:s(v)} we obtain
\begin{equation} \label{Poisson:W(z^2)}
\frac{W(z)}{z}=\frac{2}{\pi}\int_0^{\pi/2}\frac{v\tan v}{z^2+s^2(v)}\frac{ds}{dv}dv \ .
\end{equation}
Similar to \eqref{eq:v'(t)} one can find
\begin{equation} \label{eq:v'(s)}
\frac{dv}{ds}=\frac{v}{s(v)\left[v^2+(1+v\tan v)^2\right]} \ .
\end{equation}
Substituting \eqref{eq:v'(s)} and \eqref{Imz:s(v)} into \eqref{Poisson:W(z^2)}, the theorem follows.
\end{Proof}
\begin{remark}
Comparison of the formula \eqref{int:z^2} with the representation \eqref{IntegralW/z} (taken in the
right $z$-half-plane)
shows that the integrand in the former contains $z^2$ rather than $z$, which can be profitable
in using the integral representations for numerical evaluation of $W(z)$ at large $z$.
\end{remark}

\section{Poisson's integrals}
\begin{Th} \label{th:Poisson}
The following two formulae of
Poisson\footnote{The second formula is explicitly given in \cite[sec.\:80, p.\:501]{Poisson} in
terms of the tree function $T(x)$ (see \eqref{def:T(x)} and \eqref{eq:W and T} below)
and proved using the Lagrange Inversion Theorem \cite[p.\:133]{Whittaker&Watson} and
a series expansion of the logarithmic function $-\ln (1-e^{ix}\phi)$ in powers of $e^{ix}$ where
the expansion coefficients $\phi^n/n$ are exactly the coefficients of the complex exponential
Fourier series for the same function.
On the other hand, today it is well known \cite[p.\:143-145]{Caratheodory} that there is a tight
connection between the classical Poisson Formula and the Cauchy Integral Formula.
Our proof is based on the latter and thereby differs from that given in the original.}
 \mbox{\textup{\cite{Poisson}}} hold for $x\in (-1/e,e)$
\begin{equation} \label{eq:Poisson1}
W(x)= \frac{2}{\pi}\int_0^\pi \frac{\cos \sbody32\theta-xe^{-\cos\theta}\cos\left(\sbody52\theta+\sin\theta\right)}
{1-2xe^{-\cos\theta} \cos(\theta+\sin\theta)+x^2e^{-2\cos\theta}}\cos \sbody12\theta \ \mathrm{d}\theta
\end{equation}\\
\begin{equation} \label{eq:Poisson2}
W(x)=-\frac{2}{\pi}\int_0^\pi \frac{\sin \sbody32 \theta+xe^{\cos\theta}\sin\left(\sbody52\theta-\sin\theta\right)}
{1+2xe^{\cos\theta} \cos(\theta-\sin\theta)+x^2e^{2\cos\theta}}\sin\sbody12\theta\ \mathrm{d}\theta
\end{equation}
\end{Th}
\begin{Proof}
We consider the defining equation \eqref{eq:def}
for given real $z=x$
\begin{equation} \label{def:W(x)}
We^W=x
\end{equation}
and interpret it as an equation with respect to $W$.
Then we can write the equation in the form $F(W)=0$ where
\begin{equation} \label{eq:implicit}
F(\zeta)=\zeta-xe^{-\zeta} \ .
\end{equation}

Let $\Gamma$ be the positively-oriented circumference of the unit circle $\left|\zeta\right|=1$ in
the complex $\zeta$-plane and domain $G$ be the interior of $\Gamma$.
The function $F(\zeta)$ is holomorphic in $G$ and by Rouch\'e's theorem it has a single isolated
zero there when $\left|x\right|<1/e$ because in this case
$\left|-xe^{-\zeta}\right|<\left|\zeta\right|$ on $\Gamma$.
Therefore, using Cauchy's integral formula with taking $\Gamma$ for the integration contour
we can write
\begin{equation} \label{Cauchy}
W=\frac{1}{2\pi i} \int_\Gamma\frac{F^\prime(\zeta)}{F(\zeta)}\zeta d\zeta \
\end{equation}
for $\left|x\right|<1/e$.\\
Since $F^\prime(\zeta)=1+xe^{-\zeta}=1+\zeta$  by \eqref{eq:implicit} and \eqref{def:W(x)}, we obtain
\begin{equation} \label{eq:Cauchy}
W=\frac{1}{2\pi} \int_{-\pi}^\pi
\frac{e^{i\theta}(1+e^{i\theta})}{1-xe^{-\cos\theta-i(\theta+\sin\theta)}}d\theta \ ,
\end{equation}
where we set $\zeta=e^{i\theta}, -\pi\leq\theta\leq\pi$.
Separating the real and imaginary parts of the integrand in \eqref{eq:Cauchy} we find that
the former is an even function of $\theta$ whereas the latter is an odd one.
Thus, the integral of the imaginary part vanishes, as should be, and the integral of the real part
gives double the value of the integral on $[0,\pi]$.
As a result, after some arrangements, we come to integral \eqref{eq:Poisson1}.

If, instead of \eqref{def:W(x)},
we consider the equation defining the (Cayley) `tree function' $T(x)$ \cite[p.127-128]{Flajolet}
\begin{equation} \label{def:T(x)}
Te^{-T}=x
\end{equation}
and introduce function $H(\zeta)=\zeta-xe^\zeta$ in a similar way as function \eqref{eq:implicit} then
after analogous considerations and taking into account a relation
\begin{equation} \label{eq:W and T}
W(x)=-T(-x)
\end{equation}
in the final result we obtain formula \eqref{eq:Poisson2}.

Now we discuss the domain of validity of the integrals \eqref{eq:Poisson1} and \eqref{eq:Poisson2}
which is actually wider than the interval $-1/e<x<1/e$ arisen above in applying Rouch\'e's theorem.
It immediately follows from the fact that $W$ is a single valued function and therefore $F(\zeta)$,
the denominator in \eqref{Cauchy}, has a single zero in $G$ for each such $x$
that $\left|\zeta\right|<1$, i.e. for $-1/e<x<e$.
Since Rouch\'e's theorem is a consequence of the argument principle (see e.g. \cite{Markushevich}),
it would be instructive to obtain this result using the latter.
To do this, say for integral \eqref{eq:Poisson1}, we apply the argument principle to
function \eqref{eq:implicit} in case when $x>0$.
It is easy to see that function $\eta=F(\zeta)$ performs a conformal mapping of the
strip $\left\{-\infty<\Re\zeta<\infty,-\pi<\Im\zeta<\pi\right\}$, containing the entire domain $G$,
to the complex $\eta$-plane cut along two semi-infinite lines on
which $\eta=\xi\pm i\pi, \xi\geq 1+\ln x$.
We also cut the $\eta$-plane along the negative real axis to take $\left|\arg \eta\right|\leq\pi$ in
the cut plane and consider an image of $\Gamma$ which is defined by equations
\begin{subequations} \label{eq:image}
\begin{align}
\rho\cos\varphi &= \cos\theta-xe^{-\cos\theta}\cos(\sin\theta) \label{eq:image1} \ ,\\
\rho\sin\varphi &= \sin\theta+xe^{-\cos\theta}\sin(\sin\theta) \label{eq:image2} \ ,
\end{align}
\end{subequations}
where $\rho=\left|\eta\right|$ and $\varphi=\arg\eta$.

The equations \eqref{eq:image} are invariant under transformation
$\theta\rightarrow -\theta, \varphi\rightarrow -\varphi$ and describe a closed curve $\tilde{\Gamma}$
that is symmetric with respect to the real axis in the $\eta$-plane.
Suppose that while a variable point $\zeta$ moves along $\Gamma$ once in the $\zeta$-plane,
the image point $\eta=F(\zeta)$ moves on $\tilde{\Gamma}$ once in the $\eta$-plane,
making one cycle about the origin.
Then the change in argument of $\eta$ is $2\pi$ and therefore, by the argument principle
the function $F(\zeta)$ has a single zero in $G$ \cite[p.\:48]{Markushevich}.
For this it is necessary that two points on $\tilde{\Gamma}$ corresponding
to $\varphi=0$ and $\varphi=\pi$ are located on the real axis on the opposite sides of the
origin, i.e. with positive $\rho$ to be measured on the opposite rays.
Substituting $\theta=\pi$ in \eqref{eq:image} gives $\rho\cos\varphi=-1-xe$ and $\rho\sin\varphi=0$.
It can be $\rho>0$ only when $\varphi=\pi$; then $\rho=1+xe$ is positive for any $x>0$.
When $\theta=0$, we have $\rho\cos\varphi=1-x/e$ and $\rho\sin\varphi=0$.
Now $\varphi=0$ and $\rho=1-x/e>0$ when $x<e$.
Thus for $0<x<e$ the curve $\tilde{\Gamma}$ encloses the origin.
Since for these $x$ the right-hand side of equation \eqref{eq:image2} vanishes,
i.e. $\Im\eta=0$ sequentially at $\theta=-\pi, \theta=0 \mbox{ and }\theta=\pi$ as $\theta$
continuously changes from $-\pi$ to $\pi$, the curve $\tilde{\Gamma}$ is traversed once with
exactly one cycle about the origin being made.
This corresponds to the fact that the inverse of the mapping $\eta=F(\zeta)$ is continuous
in the domain bounded by the curve $\tilde{\Gamma}$ and on $\tilde{\Gamma}$ itself
and hence $\tilde{\Gamma}$ consists only of simple points \cite[Theorem 2.22]{Markushevich3}.
Thus, by the argument principle the function $F(\zeta)$ has a single zero in $G$.
Summarizing up the obtained results we conclude that the integral \eqref{eq:Poisson1} is valid
for $x\in (-1/e,e)$.
The integral \eqref{eq:Poisson2} can be considered in a similar manner.
\end{Proof}

\begin{remark}
The integral representations \eqref{eq:Poisson1} and \eqref{eq:Poisson2} can be immediately applied
to the tree function using relation \eqref{eq:W and T}.
\end{remark}
\begin{remark}
We can apply the above approach to the equation \eqref{eq:Wright}.
To eliminate a singularity at the origin we compose the integration contour of a small circle of radius,
say $r$, and the unit circle, both centered at the origin and  connected through the cut along the
negative real axis.
Then, making $r$ go to zero, we find for $0<x<e$
\begin{equation*} \label {Poisson:Wright}
W(x)=\psi(x)+
\frac{2}{\pi}\int_0^\pi \frac {\cos\displaystyle\frac{\theta}{2}+
\theta\sin\frac{3}{2}\theta-\cos\frac{3}{2}\theta\ln x}{1+2\theta\sin\theta+\theta^2-2\cos\theta\ln x +
\ln^2x} \cos\frac{\theta}{2}d\theta     \ ,
\end{equation*}
where
\begin{equation*} \label{eq:psi(x)}
\psi(x)=\int _{0}^{1}\!\frac {t-1}{\pi^2+(\ln x + t - \ln t)^2}dt \ .
\end{equation*}\\
\end{remark}

\section{Burniston-Siewert representations}
One of the analytic methods for solving transcendental equations is based on a canonical solution of
a suitably posed Riemann-Hilbert boundary-value problem \cite[p.\:183-193]{Henrici3}.
This method was found and developed by Burniston and Siewert \cite{Burniston}, its versions,
variations and applications were also considered by other authors.
By the method, a solution of a transcendental equation is represented as a closed-form integral
formula that can be regarded as an integral representation of the unknown variable.
Below we consider such integrals for the $W$ function which are based on the results of application of
the Burniston-Siewert method to solving equation \eqref{def:W(x)} obtained in
paper \cite{Anastasselou1} and the classical work \cite{Siewert}.

We start with two formulas derived in \cite{Anastasselou1} and apply them to function \eqref{eq:implicit}
\begin{equation} \label{eq:start}
W(x) = -F(0)\exp\left\{-\frac{1}{2\pi i}\int_\Gamma\frac{\ln\left( F(\zeta)/\zeta\right)}{\zeta} d\zeta \right\} \ ,
\end{equation}
\begin{equation} \label{eq:next}
W(x) = -\frac{1}{2\pi i}\int_\Gamma \ln\left( \frac{F(\zeta)}{\zeta} \right) d\zeta ,
\end{equation}
where the integration contour $\Gamma$ is the unit circle $\left| \zeta \right| = 1$ and $x\in(-1/e,e)$.
Since $F(0)=-x$ and $W(x)/x=e^{-W(x)}$, formula \eqref{eq:start} is simplified
\begin{equation} \label{eq:simple}
W(x) = \frac{1}{2\pi i}\int_\Gamma\frac{\ln\left( F(\zeta)/\zeta\right)}{\zeta} d\zeta  \ .
\end{equation}
We set $\zeta=e^{i\theta}, -\pi\leq\theta\leq\pi$.
Then, as $F(\zeta)/\zeta=F(e^{i\theta})e^{-i\theta}=R(\theta)+iI(\theta)$, where
\begin{eqnarray*}
R(\theta)=&1-xe^{-\cos\theta}\cos(\theta+\sin\theta), \\
I(\theta)=&xe^{-\cos\theta}\sin(\theta+\sin\theta),
\end{eqnarray*}
and $d\zeta/\zeta=id\theta$,
the integral \eqref{eq:simple} is reduced to
\begin{equation} \label{eq:short}
W(x) = \frac{1}{2\pi}\int_0^\pi \ln\left( R^2(\theta)+I^2(\theta) \right) d\theta \ .
\end{equation}\\
Similarly, the integral \eqref{eq:next} can be represented in the form
\begin{equation} \label{eq:W(arctan)}
W(x) = \frac{1}{2\pi}\int_0^\pi
\left\{2\arctan(I(\theta)/R(\theta))\sin\theta
-\ln\left( R^2(\theta)+I^2(\theta) \right)\cos\theta \right\} d\theta \ ,
\end{equation}
where we have taken into account that $\arg(R(\theta)+iI(\theta))=\arctan(I(\theta)/R(\theta))$ as
$R(\theta)>0$ for $0<\theta<\pi$ and $-1/e<x<e$.
We note that the integral \eqref{eq:short} has a simpler form than \eqref{eq:W(arctan)}.
Integrals similar to the above with using a function $\Phi(\zeta)=\zeta e^\zeta-x$ in our notations
instead of \eqref{eq:implicit} in formulas \eqref{eq:start}
and \eqref{eq:next} (without simplification \eqref{eq:simple}) are given in \cite{Anastasselou1}.

Thus the integrals \eqref{eq:short} and \eqref{eq:W(arctan)} representing
the principal branch of the Lambert $W$ function are valid in the domain that
contains interval $(-1/e,0)$.
However, there is one more branch that is also a real-valued function on this interval,
this is the branch $-1$ with the range $(-\infty,-1)$ (recall $W_0>-1$ and
$W_0(-1/e)=W_{-1}(-1/e)=-1$) \cite{Big paper}.
A representation of this branch can be obtained on the basis of a simple interpretation of
formula \eqref{eq:next} given in \cite{Anastasselou2}
\begin{equation} \label{eq:branch-1}
W_{-1}(x) = 1-2c-\frac{1}{2\pi i}\int_C \ln\left( \frac{F(\zeta)}{\zeta} \right) d\zeta ,
\end{equation}
where the circle $C$ is defined by equation $\left| \zeta+c \right| = c-1$
with arbitrary constant $c>1$ and $-1/e<x<-(2c-1)e^{1-2c}$.
Transformations of \eqref{eq:branch-1} leading to a definite integral
are similar to those used above to obtain the integrals \eqref{eq:short} and \eqref{eq:W(arctan)}
and skipped here together with a bulky result.

We return to the principal branch and use the result in \cite[formula(13)]{Siewert} to write \cite{Wolfram}
\begin{equation} \label{int1}
W(z) = 1+(\ln z-1)\exp\left(\frac{i}{2\pi}\int_0^\infty
\ln\left(\frac{\ln z + t - \ln t + i\pi}{\ln z + t - \ln t - i\pi}\right)\frac{dt}{1+t}\right)
\end{equation}
or

\begin{equation} \label{int:arg}
W(z) = 1+(\ln z-1)\exp\left\{ -\frac{1}{\pi}\int_0^\infty
\frac{\arg(\ln z+t-\ln t+i\pi)}{1+t}dt \right\} \ ,
\end{equation}
\\where $z\notin [-1/e,0]$.
In case of real $z=x>1/e$, when the expression $\ln z + t - \ln t$ is real and positive (for $t\in(0,\infty)$),
the formula \eqref{int:arg} is simplified and reduced to
\begin{equation} \label{int2}
W(x) = 1+(\ln x-1)\exp\left\{-\frac{1}{\pi}\int_0^\infty
\arctan\left(\frac{\pi}{\ln x + t - \ln t }\right)\frac{dt}{1+t}\right\}
\end{equation}
or, after integrating by parts
\begin{equation} \label{int3}
W(x) = 1+(\ln x-1)\exp\left\{-\int_0^\infty\frac{t-1}{\pi^2+(\ln x + t - \ln t)^2}
\frac{\ln(1+t)}{t}dt\right\} \ .
\end{equation}
\\
We emphasize that the domain $x>1/e$ of validity of the formulae \eqref{int2} and
 \eqref{int3} is different from that of \eqref{eq:short} and \eqref{eq:W(arctan)}.

For the case $x\in(-1/e,0)$, we refer the reader to \cite[formula (32)]{Siewert}
where the principal branch $W_0$ and the branch $W_{-1}$
are represented in the form of a combination of two expressions similar to
the right-hand side of \eqref{int:arg}.

\begin{remark}
We can regard the integral in the formula \eqref{int1} as an improper integral depending
on a parameter $p=\ln z$ and consider it in the limit $p\rightarrow\infty$ (when $z\rightarrow\infty$).
Since the integrand is a continuous function of two variables $t$ and $p$ in the domain under
consideration and the integral is uniformly convergent with respect to $p$,
we can take the limit under the integral sign and find that the integral vanishes as the
integrand goes to zero.
Then the formula \eqref{int1} reproduces the asymptotic result \eqref{eq:branch}.
\end{remark}
Finally we note that by use of elementary complex analysis in \cite{Kheyfits}
there is obtained a common closed form representation for all
the branches $W_k(z)$ in the complex $z$-plane through simple quadratures.

\section{Concluding remarks}
We derived various integral representations of the principal branch of the Lambert $W$ function
using different approaches.
The most part of them is associated with functions containing $W$ belonging
to various classes of functions admitting certain integral representations.
Among other classes we considered in detail the classes of Stieltjes functions and
complete monotonic functions
and by the example of functions based on $W$ in fact demonstrated different ways
to establish belonging of a function to these classes.

Besides their own importance the derived integral representations have some applications.
One of them has been mentioned in connection with finding Nuttall-Bouwkamp integral \eqref{eq:Nutend}.
Other definite integrals appear in taking the obtained integrals with a particular value of $z$.
For example, integrals \eqref{IntegralW/z}, \eqref{eq:1/1+W}, \eqref{eq:1/W}, \eqref{int:z^2} taken
at $z=e$ yield respectively
\begin{equation*}
\int_0^\pi\frac{v^2+(1-v\cot v)^2}{1+v\csc(v)e^{-(1+v\cot v)}}dv=\pi \ ,
\end{equation*}
\\
\begin{equation*}
\int_0^\pi\frac{v^2+(1-v\cot v)^2}{v\csc(v)\left(v\csc(v)+e^{1+v\cot v}\right)}dv =
\frac{e-1}{e}\,\pi \ ,
\end{equation*}
\begin{equation*}
\int_0^\pi\frac{dv}{1+e^{1+v\cot v}\sin v/v}=\frac{1}{2}\,\pi \ ,
\end{equation*}
\begin{equation*}
\int_0^{\pi/2}\frac{\left[v^2+(1+v\tan v)^2\right]v\sec(v)e^{v\tan v-1} }
{1+v^2\sec^2(v)e^{2(v\tan v-1)}}\tan v\,dv=\frac{1}{2}\,\pi \ .
\end{equation*}
Another advantage that can be taken of the obtained results is based on
a comparison between different representations of the same function.
This reveals equivalent forms of the involved integrals.
In addition, since some of the integrals are simpler than others,
such equations can be regarded as a simplification of the latter.
For example, equating integrals \eqref{Pick3} and \eqref{IntegralW/z} shows that the former
can be simplified and reduced to the latter.

At last we mention that the Pick representations \eqref{Pick1}, \eqref{Pick2},
\eqref{Pick3}, and \eqref{Pick4}
can be considered as integrals expressing properties of the kernel $K(z,v)$ defined by \eqref{K(z,v)}.


\end{document}